\documentclass[letterpaper, 10 pt, conference]{ieeeconf}  % Comment this line out
\usepackage{graphicx}
\usepackage{amsmath}
\usepackage{amsfonts}                                                     % paper
\usepackage{amssymb}

\IEEEoverridecommandlockouts                              % This command is only
\title{\LARGE \bf
Identification of Structured LTI MIMO State-Space Models
}

\author{Chengpu Yu,  Michel Verhaegen, Shahar Kovalsky and Ronen Basri % <-this % stops a space
\thanks{C. Yu and M. Verhaegen are with the Delft Center for Systems and Control, Delft University, Delft 2628CD, Netherlands (c.yu-4@tudelft.nl, m.verhaegen@tudelft.nl)}
\thanks{S. Kovalsky and R. Basri are with Weizmann Institute of Science, Department of Computer Science and Applied Mathematics (shahar.kovalsky@weizmann.ac.il, ronen.basri@weizmann.ac.il) }
 \thanks{The work of C. Yu and M. Verhaegen is sponsored by the European Research Council, Advanced Grant Agreement No. 339681.}
}
\begin{document}

\maketitle
\thispagestyle{empty}
\pagestyle{empty}

%%%%%%%%%%%%%%%%%%%%%%%%%%%%%%%%%%%%%%%%%%%%%%%%%%%%%%%%%%%%%%%%%%%%%%%%%%%%%%%%
\begin{abstract}

The identification of structured state-space model has been intensively studied for a long time but still has not been adequately addressed. The main challenge is that the involved estimation problem is a non-convex (or bilinear) optimization problem. This paper is devoted to developing an identification method which aims to find the global optimal solution under mild computational burden. Key to the developed identification algorithm is to transform a bilinear estimation to a rank constrained optimization problem and further a difference of convex programming (DCP) problem. The initial condition for the DCP problem is obtained by solving its convex part of the optimization problem which happens to be a nuclear norm regularized optimization problem. Since the nuclear norm regularized optimization is the closest convex form of the low-rank constrained estimation problem, the obtained initial condition is always of high quality which provides the DCP problem a good starting point. The DCP problem is then solved by the sequential convex programming method. Finally, numerical examples are included to show the effectiveness of the developed identification algorithm.

\end{abstract}

\section{Introduction}
This paper studies the identification of structured state-space systems. In the literature, there are two kinds of parameterized state-space models \cite{ljung1998system}. One is the black-box model for which the associated system matrices are fully parameterized. The other is the gray-box model (also called structured model in this paper), where only parts of system matrices are parameterized as non-zero parameters. For the identification of black-box LTI system models, the subspace identification methods are commonly adopted as they generally yield reliable identification results \cite{verhaegen2007filtering}. Other identification methods for black-box systems, such as the regularized Gauss-Newton method and the gradient projection method, are also widely applied in practice \cite[Chapter 7]{verhaegen2007filtering}. Although the gray-box model generally contains less parameters than the black-box model, it is usually difficult to obtain reliable identification results. The main reason is that the associated identification problem for the gray-box model is always non-convex and it may have many local optimal points. In this regard, it is difficult to find a global optimal solution under mild computational burden. Generally, the gray-box model is solved by either the output-error method or the prediction-error method for which the involved nonlinear optimization problem is always tackled by gradient-based algorithms \cite{verhaegen2007filtering,lyzell2012structural}.

For the structured system model, apart from the direct identification method which is to estimate parameters directly from the system inputs and outputs, the indirect identification framework has also been intensively investigated recently, e.g., see \cite{lyzell2009handling,prot2012null,mercere2014identification}. This indirect identification framework contains two steps. First, the system matrices of the concerned state-space model are identified up to a similarity transformation using classical subspace identification methods, such as MOESP and N4SID \cite{verhaegen2007filtering}. Second, the system parameters are determined according to the system parametrization pattern and the obtained system matrices in the first step.

The second step of the indirect identification for a structured state-space system turns out to be a nonlinear estimation problem which has a smaller scale than that in the direct identification method. More specifically, the involved nonlinear estimation problem is bilinear, i.e. the system parameters and the similarity transformation matrix are coupled together. To solve such a kind of bilinear estimation problem, the alternating minimization algorithm is adopted in \cite{xie2002estimate}, which tries to minimize the objective function with respect to either the model parameters or the similarity transformation matrix. The main shortage of the alternating minimization algorithm is that it is likely to get stuck in local optimal points under randomly generated initial conditions. To cope with the initialization problem, the bilinear estimation problem is reformulated as a sum-of-squares which is then solved by semi-definite programming \cite{ljung2003initialization}. Due to the complexity of the sum-of-squares decomposition, this method can only be applied to a small-scale structured system. Since many variables are involved in this identification method, the estimate of the similarity transformation tends to be ill-conditioned.

Recently, in order to robustly solve the bilinear estimation problem for the identification of a structured system, a null-space-based technique is provided in \cite{prot2012null}. It first stacks all the variables together, no matter independent or coupled, and then compute a subspace which the augmented variable vector lies in. By representing all the variables using the linear combination of the subspace basis, the dimension of the associated non-convex optimization problem can be reduced, hence decreasing computational burden and improving identification accuracy. In practical simulations, the estimate of the similarity transformation matrix may be ill conditioned or even singular when the adopted algorithm get stuck to local (not global) optimal points. To overcome this problem, the condition number of the similarity transformation matrix is constrained when dealing with the bilinear estimation problem \cite{mercere2014identification}, and a combination of BFGS method and spectral bundle algorithm is implemented to solve such a highly nonlinear optimization problem.

In this paper, we present a new identification method, which is a combination of the over-parametrization technique and low-rank constrained optimization. The main idea of this method is that the considered nonlinear estimation problem can be equivalently transformed to a linear estimation problem with a low rank constraint. The relaxation from a nonlinear estimation problem to a linear estimation problem introduces a number of redundant variables, while the low rank constraint reduces this redundancy. The low-rank constrained optimization problem is generally NP hard. However, for the investigated identification problem, the associated rank constrained optimization problem can be casted to a difference of convex programming (DCP) problem, which can be solved using many existing methods \cite{boyd2008sequential,tao2005dc}. The initial point for solving the DCP problem can be obtained by solving the convex part of the DCP problem. Interestingly, the convex part of the DCP problem is exactly a nuclear norm regularized optimization problem, which can be commonly employed as a heuristic way to solve low-rank constrained optimization problems. By solving the nuclear-norm regularized optimization problem, a good candidate of the initial point can usually be found. Then, by iteratively solving the DCP problem, the identification performance can be gradually improved, which usually yields a global optimal solution in practice. Compared with the methods in \cite{prot2012null,ljung2003initialization}, the reduction of the variable redundancy using the low-rank constraint is the key contribution of the present paper.

The rest of the paper is organized as follows. Section \ref{sec1} provides the problem formulation and a discussion on system identifiability. Section \ref{sec2} presents a new identification method for gray-box models. Section \ref{sec3} gives simulation examples to validate the proposed identification method. Finally, concluding remarks are given in Section \ref{sec4}.

\section{Problem statement and system identifiability}\label{sec1}
In this paper, we consider the identification of parameterized LTI state-space system model. Let $\theta\in\mathbb R^{q}$ be a parameter vector. The concerned state-space is described as
\begin{equation}\label{equ1}
  \begin{split}
    &x(k+1)=A(\theta)x(k)+B(\theta)u(k)\\
    &y(k)=C(\theta)x(k)+w(k),
  \end{split}
\end{equation}
where $u(k)\in\mathbb R^{m},x(k)\in\mathbb R^{n}, y(k)\in\mathbb R^{p}$ and $w(k)\in\mathbb R^{p}$ are system input, state, output and measurement noise, respectively. The system matrices $A(\theta),B(\theta)$ and $C(\theta)$ are assumed to be affine functions with respect to $\theta=\left[\theta_1,\cdots,\theta_q\right]^T$:
\begin{equation}
  \begin{split}
    &A(\theta)=A_0+A_1\theta_1+\cdots+A_q\theta_q,\\
    &B(\theta)=B_0+B_1\theta_1+\cdots+B_q\theta_q,\\
    &C(\theta)=C_0+C_1\theta_1+\cdots+C_q\theta_q.
  \end{split}
\end{equation}
In practice, the parameter vector $\theta$ in \eqref{equ1} often has physical interpretation, such as the Newton's or Kirchoff's laws \cite[Example 7.1]{verhaegen2007filtering} and structure of a compartmental model \cite{bellman1970structural}. Thus, the estimation of these parameters is practically meaningful.

For the system model in \eqref{equ1}, the system input is assumed to be persistently exciting and the system model is minimal (observable and controllable). When the measurement noise $w(k)$ is uncorrelated with the system input $u(k)$, the associated system matrices can be consistently estimated up to a similarity transformation using the subspace identification method \cite{ljung1998system,verhaegen2007filtering}. Suppose that we have obtained a correct estimate of system matrices, up to a non-singular similarity transformation, denoted as $\hat A,\hat B$ and $\hat C$. Then, the following equation holds
\begin{equation}\label{equ2}
  \begin{split}
    &\hat A=T^{-1}A(\theta)T\\
        &\hat B=T^{-1}B(\theta)\\
        &\hat C=C(\theta)T,
  \end{split}
\end{equation}
where $T\in\mathbb R^{n\times n}$ is a non-singular ambiguity matrix. In the above equation, both $\theta$ and $T$ are unknown.  Before proceeding to solve the above nonlinear equation, we shall discuss the system identifiability of \eqref{equ1} or the solution uniqueness of \eqref{equ2}.

In order to check the identifiability of the state-space model in \eqref{equ1}, a direct way is to provide the corresponding system transfer function in terms of $\theta$ and check whether the involved parameters can be uniquely determined from the system impulse response \cite{bellman1970structural,norton2010identification}. There are generically two factors affecting the system identifiability: the number of involved parameters and the parametrization pattern of the model. For the state-space model in \eqref{equ1}, if there are more than $n(p+m)$ parameters, the system is unidentifiable \cite{ljung1998system,verhaegen2007filtering}. Thus, a necessary condition is that the number of parameters should be less than $n(p+m)$. As an extreme case, if all the system matrices are fully parameterized, then they are unidentifiable. However, less than $n(p+m)$ does not necessarily result in the identifiability of the system model \eqref{equ1}. In fact, the identifiability also depends on the structure of the concerned system model. It is shown in \cite{bellman1970structural} that the state-space model described by the following parameterized system matrices is unidentifiable even though the number of parameters is strictly less than $n(m+p)$:
\begin{equation}
  \begin{split}
    &A(\theta)=\left[\begin{array}{ccc}
      -\theta_1 &\theta_3&0\\
      \theta_1&-(\theta_2+\theta_3)&\theta_4\\
      0&\theta_3&-\theta_4
    \end{array}\right], B(\theta)=\left[\begin{array}{c}
      0\\
      0\\
      1
    \end{array}\right],\\
    & C(\theta)=\left[\begin{array}{ccc}
      1 &0&0
    \end{array}\right].\nonumber
  \end{split}
\end{equation}
It has been shown in \cite{bellman1970structural} that the following system structures are identifiable: the diagonal structure, companion structure, cascaded compartmental structure with only system input and output taking place at the last compartment. It is also straightforward that any system model derived from the above identifiable structure under some similarity transformation is also identifiable.

Since this paper mainly focuses on developing a robust identification method for structured state-space system models, the model structure in \eqref{equ1} is assumed to be identifiable. The above assumption indicates that the parameter vector $\theta$ can be uniquely determined from \eqref{equ2}, as well as the ambiguity matrix $T$. In addition, by the assumption that the concerned system is minimal, we can see that $T$ is a nonsingular matrix. Although $\theta$ and $T$ are uniquely determined by \eqref{equ2}, it is usually difficult to find the exact solution of the nonlinear estimation problem.

Equation \eqref{equ2} can be reformulated as
\begin{equation}\label{equ3}
  \begin{split}
    &T\hat A=A(\theta)T\\
        &T\hat B=B(\theta)\\
        &\hat C=C(\theta)T.
  \end{split}
\end{equation}
A basic approach to estimate $T$ and $\theta$ is to solve the following optimization problem \cite{xie2002estimate}:
\begin{equation}\label{equ4}
  \begin{split}
    &\arg\min_{\theta, T}\quad h(\theta,T)\\
    &h(\theta,T)=\|T\hat A-A(\theta)T\|_F^2+\|T\hat B-B(\theta)\|_F^2\\
    &\quad\quad +\|\hat C-C(\theta)T\|_F^2.
  \end{split}
\end{equation}

 Since the above estimation problem is non-convex, it is difficult to obtain a global optimal solution under mild computational burden. In order to obtain a solution close to the global optimal point, selection of the initial condition is crucial. Although a semi-definite programming method is provided in \cite{ljung2003initialization} to find a proper initial point, the involved the sums-of-squares decomposition of \eqref{equ4} is nontrivial and computationally expensive. Also, an increased number of decision variables are involved, which makes the estimation problem more under-determined. In this paper, we shall develop an effective and robust identification algorithm for structured state-space models. The associated algorithm relies on the combination of the over-parametrization technique and the low-rank constraint. The over-parametrization technique suggests to introduce more variables to transform a non-linear estimation problem into a linear estimation problem, as done in \cite{ljung2003initialization}; however, the low-rank constraint is employed to exploit the relations among all involved variables and to reduce their associated redundancy.

\section{Gray-box system identification}\label{sec2}
In this subsection, we aim to develop an algorithm to estimate $\theta$ and $T$ from equation \eqref{equ3}. Denote by $\otimes$ and $\textrm{vec}(\cdot)$ respectively the Kronecker product and the vectorization operators. The vectorized form of equation \eqref{equ3} is shown as follows:
  \begin{equation}\label{equ5}
    \underbrace{\left[\begin{array}{c}
      \hat A^T\otimes I-I\otimes A(\theta)\\
      \hat B^T\otimes I\\
      I\otimes C(\theta)
    \end{array}\right]}_{M(\theta)}\textrm{vec}(T)=\underbrace{\left[\begin{array}{c}
      0\\
      \textrm{vec}(B(\theta))\\
      \textrm{vec}(\hat C)
    \end{array}\right]}_{N(\theta)},
  \end{equation}
  where $M(\theta)\in\mathbb R^{(n^2+mn+pn)\times n^2}$ and $N(\theta)\in\mathbb R^{(n^2+mn+pn)}$ are introduced for notational simplicity. Since $A(\theta),B(\theta)$ and $C(\theta)$ are affine functions in terms of $\theta$, $M(\theta)$ and $N(\theta)$ in the above equation are affine as well. Denote the affine expressions of $M(\theta)$ and $N(\theta)$ as follows:
\begin{equation}
    \begin{split}
      &M(\theta)=M_0+M_1\theta_1+\cdots+M_q\theta_q,\\
      &N(\theta)=N_0+N_1\theta_1+\cdots+N_q\theta_q,
    \end{split}
\end{equation}
where $\{M_i\}_{i=0}^q$ and $\{N_i\}_{i=0}^q$ are constant coefficient matrices of $M(\theta)$ and $N(\theta)$, respectively. Let $\tau=\textrm{vec}(T)$. Then equation \eqref{equ5} can be equivalently expressed as
\begin{equation}\label{equ6}
  M_0\tau+M_1\tau\theta_1+\cdots+M_q\tau\theta_q=N_0+N_1\theta_1+\cdots+N_q\theta_q.
\end{equation}
From the above equation, we can find that the associated estimation problem is bilinear. More specifically, the variables of $\theta$ and $\tau$ are coupled together. To cope with this, inspired by the over-parametrization technique, we introduce the following redundant variables:
\begin{equation}\label{equ8}
    \begin{split}
      &\vartheta_1=\tau\theta_1,\\
      &\quad\quad\vdots\\
      &\vartheta_q=\tau\theta_q,
    \end{split}
\end{equation}
where $\vartheta_i\in\mathbb R^{n^2}$ for $i=1,\cdots,q$. Then equation \eqref{equ6} is formed as
\begin{equation}\label{equ7}
  M_0\tau+M_1\vartheta_1+\cdots+M_q\vartheta_q=N_0+N_1\theta_1+\cdots+N_q\theta_q,
\end{equation}
where $\{\vartheta_i\}_{i=1}^q$ and $\{\theta_i\}_{i=1}^q$ are unknown variables. The above equation is linear with respect to the unknown variables, but there are more variables than equations; thus, the associated estimation problem is under-determined. Constraining or regularizing the variables are typically used to mitigate this under-determinedness. Here, by exploring the inherent relations among the variables in \eqref{equ8}, we can find that the following composed matrix is a rank one matrix:
\begin{equation}
  H(\vartheta,\theta)=\left[\begin{array}{cccc}
    \tau&\vartheta_1&\cdots&\vartheta_q\\
    1&\theta_1&\cdots&\theta_q
  \end{array}\right],
\end{equation}
where $\vartheta=\left[\begin{array}{ccc}
\vartheta_0&\cdots&\vartheta_q
\end{array}\right]\in\mathbb R^{n^2\times q}$.

As will be shown in the following proposition, equation \eqref{equ6} is equivalent to equation \eqref{equ7} with a rank constraint.
\newtheorem{pro1}{Proposition}
\begin{pro1}\label{proposition1}
  Assume  that there exists a unique solution pair $(\theta, \tau)$ for equation \eqref{equ6}. The variables $(\vartheta,\theta)$ can be uniquely determined from equation \eqref{equ7} subject to the following rank constraint:
  \begin{equation}
    \textrm{rank}\left[H(\vartheta,\theta)\right]=\textrm{rank}\left[\begin{array}{cccc}
    \tau&\vartheta_1&\cdots&\vartheta_q\\
    1&\theta_1&\cdots&\theta_q
  \end{array}\right]=1.
  \end{equation}
  Furthermore, $(\theta, \tau)$ can be uniquely determined from $(\vartheta,\theta)$.
\end{pro1}

\begin{proof}
  To prove the above proposition, it suffices to prove that the variables $(\theta,\tau)$ and $(\vartheta,\theta)$ are in one-to-one mapping under the rank constraint. It is obvious that the variables $(\vartheta,\theta)$ can be uniquely determined from $(\theta,\tau)$. Next, we will show that $(\theta,\tau)$ can be determined from $(\vartheta,\theta)$ under the rank constraint as well.

   Since $H(\vartheta,\theta)$ is a rank one matrix, its SVD decomposition can be represented as $H(\vartheta,\theta)=uv^T$ where $u\in\mathbb R^{(n^2+1)}$ and $v\in\mathbb R^{(q+1)}$. It follows that
  $uv^T=\left[\begin{array}{cccc}
    \tau&\vartheta_1&\cdots&\vartheta_q\\
    1&\theta_1&\cdots&\theta_q
  \end{array}\right]=\left[\begin{array}{c}
    \tau\\
    1
  \end{array}\right]\left[\begin{array}{cc}
    1&\theta^T
  \end{array}\right]$ with $\tau\in\mathbb R^{n^2}$.
  It then follows that $u=\alpha\left[\begin{array}{c}
    \tau\\
    1
  \end{array}\right]$ and $v=\frac{1}{\alpha}\left[\begin{array}{c}
    1\\
    \theta
  \end{array}\right]$ with $\alpha$ a scalar ambiguity. Obviously, $\theta$ and $\tau$ can be respectively determined from $v$ and $u$. Therefore, $(\theta,\tau)$ can be determined from $(\vartheta,\theta)$.
\end{proof}

By the above proposition, instead of dealing with the bilinear estimation of \eqref{equ6}, we shall consider the following equivalent estimation problem:
\begin{equation}\label{equ88}
    \begin{split}
      \min_{\vartheta,\theta}& \quad \|\sum_{i=0}^q M_i\vartheta_i-N_0-\sum_{i=1}^q N_i\theta_i\|_2^2\\
      s.t.&\quad \textrm{rank}\left[H(\vartheta,\theta)\right]=1.
    \end{split}
\end{equation}

The above rank constrained optimization problem is NP hard, which is difficult to find a global optimal solution under mild computational burden. One heuristic way for the above optimization problem is to solve the following nuclear norm regularized optimization problem:
\begin{equation}\label{equ9}
  \min_{\vartheta,\theta}\quad \|\sum_{i=0}^q M_i\vartheta_i-N_0-\sum_{i=1}^q N_i\theta_i\|_2^2+\lambda\|H(\vartheta,\theta)\|_*,
\end{equation}
where $\lambda$ is a regularization parameter. Since the above optimization problem is convex, we can reliably obtain an approximate estimate of $(\vartheta,\theta)$.

Denote by $\sigma_i\left(H(\vartheta,\theta)\right)$ the $i$-th largest singular value of $H(\vartheta,\theta)$ for $i=1,\cdots, q+1$. Let
\[
f_{\kappa}\left(H(\vartheta,\theta)\right)=\sum_{i=1}^{\kappa}\sigma_i\left(H(\vartheta,\theta)\right) \textrm{ for }\kappa=1,\cdots,q+1.
 \]
It is remarked that $f_{\kappa}\left(\cdot\right)$ is a Ky Fan $\kappa$-norm \cite{bhatia2013matrix}. Inspired by the truncated nuclear norm method in \cite{hu2013fast}, the rank constraint in \eqref{equ88} can be replaced by the following constraint:
\[
\sum_{i=2}^{q+1}\sigma_i\left(H(\vartheta,\theta)\right)=f_{q+1}\left(H(\vartheta,\theta)\right)-f_1\left(H(\vartheta,\theta)\right)=0.
  \]
As a consequence, we try to solve the following unconstrained optimization problem:
\begin{equation}\label{equ10}
\begin{split}
  \min_{\vartheta,\theta}&\quad \|\sum_{i=0}^q M_i\vartheta_i-N_0-\sum_{i=1}^q N_i\theta_i\|_2^2\\
  &\quad\quad+\lambda\left[f_{q+1}\left(H(\vartheta,\theta)\right)-f_1\left(H(\vartheta,\theta)\right)\right].
\end{split}
\end{equation}
According to the definition of $f_{\kappa}(\cdot)$, we observe that
\[
\begin{split}
&f_{q+1}\left(H(\vartheta,\theta)\right)=\|H(\vartheta,\theta)\|_*\\
&f_1\left(H(\vartheta,\theta)\right)=\|H(\vartheta,\theta)\|_2.
\end{split}
\]
 Then equation \eqref{equ10} is therefore equivalent to
\begin{equation}\label{equ11}
\begin{split}
  &\min_{\vartheta,\theta}\quad \|\sum_{i=0}^q M_i\vartheta_i-N_0-\sum_{i=1}^q N_i\theta_i\|_2^2\\
  &\quad\quad\quad+\lambda\left(\|H(\vartheta,\theta)\|_*-\|H(\vartheta,\theta)\|_2\right).
\end{split}
\end{equation}
Due to the convex properties of the nuclear norm and the spectral norm, it is easy to see that the above optimization is a difference of convex programming (DCP) problem, which can be efficiently solved by sequential convex programming method \cite{boyd2008sequential} or the difference convex algorithm (DCA) \cite{tao2005dc}. Since these two algorithms have similar convergence properties \cite{lanckriet2009convergence}, we shall apply the sequential convex programming method to solve the optimization problem in \eqref{equ11}.

Denote by  $(\hat \vartheta^k,\hat\theta^k)$ the estimate of $(\vartheta,\theta)$ at the $k$-th iteration. Let $\hat u^k$ and $\hat v^k$ be the left and right singular vectors associated with the largest singular value of $H(\hat\vartheta^k,\hat\theta^k)$, respectively. Then, linearizing the concave term of the objective function in \eqref{equ11} and leaving out the constant terms yields
\begin{equation}\label{equ12}
\begin{split}
  \min_{\vartheta,\theta} &\quad\|\sum_{i=0}^q M_i\vartheta_i-N_0-\sum_{i=1}^q N_i\theta_i\|_2^2\\
  &+\lambda\left(\|H(\vartheta,\theta)\|_*-(\hat u^k)^TH(\vartheta,\theta)\hat v^k\right).
\end{split}
\end{equation}
Since \eqref{equ12} is convex, the estimate $(\hat\vartheta^{k+1},\hat\theta^{k+1})$ can be reliably obtained using existing convex optimization tools. It is noteworthy that iterative minimization of \eqref{equ12} is the sequential convex programming method.

Since the DCP is still a non-convex optimization problem, the obtained solution using the adopted sequential convex programming method depends on the initial condition. To cope with this, the initial condition is obtained by solving the nuclear norm regularized optimization problem in \eqref{equ9}. In fact, the optimization problem in \eqref{equ9} is the convex part of \eqref{equ11}. Since the nuclear norm is the convex envelope of the low-rank constraint on the unit ball of the operator norm \cite{recht2010guaranteed}, the solution to \eqref{equ9} is usually a good candidate for the starting point of the sequential convex programming method. Numerous simulation results show that based on such an initial condition, the sequential convex programming method usually leads to a global optimal solution.

By the identifiability assumptions on the model structure, the matrix $T$ has to  be nonsingular at the global optimal point. However, getting stuck in local optimal points may cause the estimate of matrix $T$ to be ill-conditioned and even singular \cite{ljung2003initialization}. Due to the continuity of the equation in \eqref{equ2}, there always exists a region around the true solution such that the global optimal solution can be obtained once the initial point of the DCP programming approach falls in this region. In other words, if the initial solution is very close the global optimal point, the estimate of $T$ will be automatically nonsingular.

\section{Numerical simulation}\label{sec3}
In this section, two examples are simulated to show the performance of the proposed identification method. The first example comes from \cite{bellman1970structural}, which is to identify compartmental structures of a networked system. The second one is to estimate a state-space represented printer belt derive model with its parameters having physical iterpretation \cite{mercere2014identification}. The system orders of these two examples are set to three.

To evaluate the performance of the developed identification method, we use the root normalized mean square error criterion which is defined as
\begin{equation}
  \textrm{rNMSE}=\sqrt{\frac{1}{K}\sum_{i=1}^K \frac{\|\hat\theta^i-\theta^*\|^2}{\|\theta^*\|^2}},
\end{equation}
where $K$ is the number of Monte-Carlo trials which is set to 50 in the following simulations, $\hat\theta^i$ and $\theta^*$ are the $i$-th estimated and the true parameters, respectively.

In the simulations, both the system input and measurement noise are white noises which are generated independently. Therefore, the persistent excitation property for the system input and the consistent estimation conditions are satisfied with probability one. To estimate the system matrices, we employ the MOESP subspace identification method \cite{verhaegen2007filtering,verhaegen2007filtering1}. In addition, the convex optimization problems involved are solved using the CVX toolbox. The regularization parameter $\lambda$ in \eqref{equ9} is empirically set to $\lambda=10^{-3}$. For each example, we show the identification performance at different noise levels. Moreover, we choose 400 input-output data for the system identification in each Monte-Carlo trial. The stopping criterion for the proposed iterative identification method is set to
\[
    \frac{\|\hat\theta^{i+1}-\hat\theta^i\|}{\|\hat\theta^i\|}\leq 10^{-6}.
\]

\textbf{Example 1:} In this example, the system model comes from a networked system with compartmental structures \cite{bellman1970structural}. The involved system matrices are parameterized as follows:
\begin{equation}
\begin{split}
  A(\theta)&=\left[\begin{array}{ccc}
    -\theta_1&\theta_3\\
    \theta_1&-\theta_2-\theta_3&\theta_4\\
    &\theta_2&-\theta_4
  \end{array}\right], B(\theta)=\left[\begin{array}{c}
    0\\
    0\\
    1
  \end{array}\right],\\
  C(\theta)&=\left[\begin{array}{ccc}
    0&0&1
  \end{array}\right].
\end{split}
\end{equation}
Compared with the unidentifiable example in Section \ref{sec1}, the matrix $C$ in this example has a different value, and the concerned system model here has been proven to be identifiable \cite{bellman1970structural}. The true value of the parameter vector is set to
\[\theta=\left[\begin{array}{cccc}
  -0.394 &  -0.893   & 0.325 &   0.383
\end{array}\right].
\]

For the ease of reference, solving the optimization problem in \eqref{equ9} to estimate the system parameters is called the \emph{\textbf{nu}clear-\textbf{n}orm (NUN) method}. The obtained estimate of the system parameters is the initial point for the  \emph{DCP method} presented in Section \ref{sec2}. For the comparison purpose, based on initial point obtained by the \emph{NUN method}, we carry out the \emph{\textbf{a}lternating \textbf{mi}nimization (AMI) method} \cite{xie2002estimate} to iteratively estimate the system parameters. For the purpose of fair comparison, both the \emph{DCP method} and the \emph{AMI method} run for 100 iterations.

\begin{figure}[htb]
  \centering
  \includegraphics[scale=0.55]{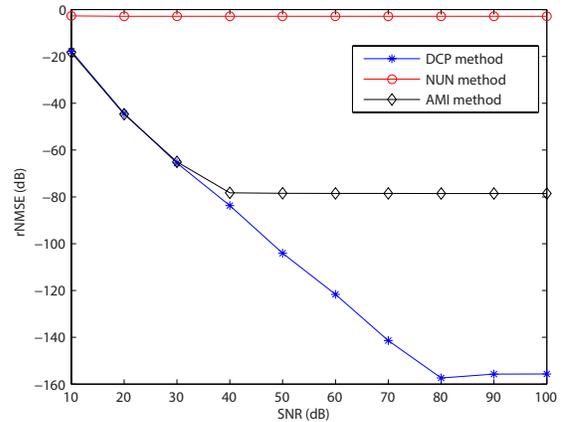}
  \caption{Example 1: identification performance against the SNR}\label{figure1}
\end{figure}
\begin{figure}[htb]
  \centering
  \includegraphics[scale=0.55]{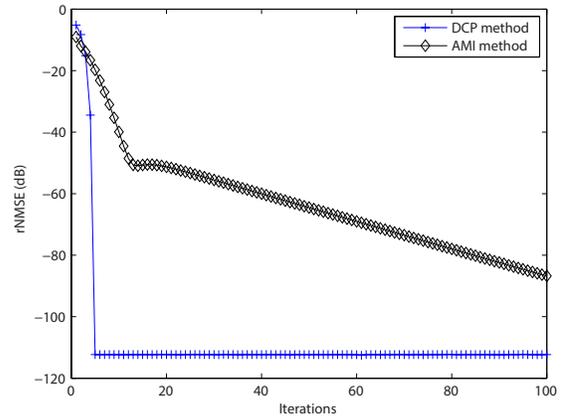}
  \caption{Example 1: identification performance against the number of iterations for a Monte-Carlo trial at SNR=50 dB.}\label{figure2}
\end{figure}

Fig. \ref{figure1} shows the identification performance of three different algorithms. We can find that the rNMSE curve corresponding to the \emph{NUN method} changes slightly along with the increase of SNR. More specifically, the estimates obtained by the \emph{NUN method} fluctuate around the point
\[
\hat{\theta}^0=\left[\begin{array}{cccc}
  -0.091& -0.304& 0.156& -0.018
\end{array}\right].
\]
 However, those rNMSE curves of the \emph{DCP method} and the \emph{AMI method} decays along with the increase of SNR. From the comparison between the \emph{DCP} and \emph{AMI methods}, we can see that the performance of the \emph{AMI method} is worse than that of the \emph{DCP method}.

Fig. \ref{figure2} shows the identification performances of the \emph{DCP } and \emph{AMI methods} on a  Monte-Carlo trial at SNR=50 dB, where we can find that the \emph{AMI method} has not met the final stopping criterion up to 100 iterations. This explains the worse performance of the \emph{AMI method} shown in Fig. \ref{figure1}. On the other hand, we can see that the \emph{DCP method} can meet the stopping criterion within 5 iterations in this simulation example, indicating that the \emph{DCP method} converges faster.

 It is remarked that the same initial point is provided for the \emph{DCP} and \emph{AMI methods} at the zeroth iteration. However, the rNMSE curves in Fig. \ref{figure2} start from the first iteration. This explains why their corresponding rNMSE values at the first iteration are slightly different.

\textbf{Example 2:} The example is extracted from \cite{mercere2014identification}, which aims to identify the physical parameters governing a printer belt drive system. The associated system matrices are parameterized as follows:
\begin{equation}
  A(\theta)=\left[\begin{array}{ccc}
    0&-1&0.15\\
    0.2&0&0\\
    \theta_1&\theta_2&\theta_3
  \end{array}\right], B=\left[\begin{array}{c}
    0\\
    0\\
    \theta_4
  \end{array}\right],C=\left[\begin{array}{c}
    0\\
  1\\
  0
  \end{array}\right]^T.
\end{equation}
The true value of $\theta$ is set to $\theta=\left[\begin{array}{cccc}
  -0.537   & 0.567 &  -0.363  & 0.156
\end{array}\right]$. All other simulation settings are the same as the previous example.

\begin{figure}[htb]
  \centering
  \includegraphics[scale=0.55]{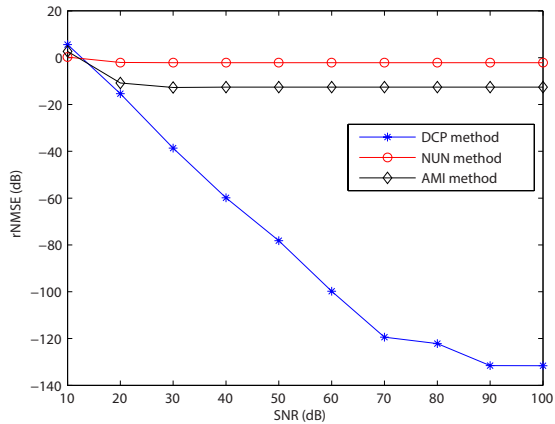}
  \caption{Example 2: identification performance against the SNR.}\label{figure3}
\end{figure}

\begin{figure}[htb]
  \centering
  \includegraphics[scale=0.55]{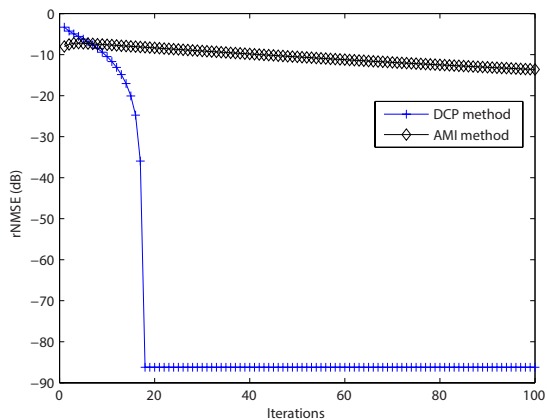}
  \caption{Example 2: identification performance against the number of iterations for a Monte-Carlo trial at SNR=50 dB.}\label{figure4}
\end{figure}

The identification performances of three different methods are shown in Fig. \ref{figure3}. The obtained estimates of the system parameters by the \emph{NUN method} fluctuate around the point
\[
\hat\theta^0=\left[\begin{array}{cccc}
 -0.143& 0.006&  -0.340& 0.157
\end{array}\right].
\]
Analogous to the previous example, the rNMSE curves of the \emph{DCP} and \emph{AMI methods} decay along with the increase of SNR. From Fig. \ref{figure4}, we can find that the \emph{DCP method} can meet the stopping criterion within 20 iterations while the \emph{AMI method} is not able to meet the stopping criterion up to 100 iterations.

\section{Conclusion}\label{sec4}
In this paper, the identification of structured state-space model has been studied and a new identification method has been devised. Key to the proposed identification method is that the concerned bilinear estimation is transformed to a rank constrained estimation problem and further a DCP problem, for which the initial condition can be robustly estimated by solving its convex optimization part. Simulation results show that the developed method usually converges to the global optimal solution. While the current approach is dedicated to the identification of structured systems, it can be straightforwardly generalized to  solve many other bilinear estimation problems, such as the output feedback controller design, blind image deconvolution, source signal separation and so on.

 \bibliographystyle{ieeetr}
 \bibliography{paper_CDC2015}

\end{document}